\documentclass[a4paper,10pt]{article}

\usepackage{authblk}

\usepackage[utf8]{inputenc} 
\usepackage[T1]{fontenc}    
\usepackage{indentfirst}
\usepackage{lmodern}

\usepackage{mathtools}
\usepackage{verbatim}
\usepackage{url}            
\usepackage{booktabs}       
\usepackage{amsfonts}       
\usepackage{amsmath}
\usepackage{amssymb}
\usepackage{amsthm}
\usepackage{microtype}      
\usepackage[a4paper]{geometry} 
\usepackage{enumitem}
\usepackage{xcolor}
\usepackage{esvect}
\usepackage{bbm}

\usepackage[linktocpage=true]{hyperref}
\hypersetup{
    colorlinks=true,
    linkcolor=blue,
    urlcolor=blue
}

\DeclareMathOperator{\Gram}{Gram}
\DeclareMathOperator{\Span}{span}

\usepackage{graphicx}
\usepackage[makeroom]{cancel}
\usepackage{tikz}

\usepackage[inkscapeformat=png]{svg}
\usepackage[normalem]{ulem}

\title{Almost sure linear independence of absolutely continuous Hilbert space-valued random vectors with respect to a special class of Hilbert space probability measures}
\author{Nizar El Idrissi* \\
nizar.elidrissi@uit.ac.ma \and Hicham Zoubeir* \\
hzoubeir2014@gmail.com}

\newcommand{\Addresses}{{
  \bigskip
  \footnotesize

  \textbf{Nizar El Idrissi.}
  \par\nopagebreak Laboratoire : Equations aux dérivées partielles, Algèbre et Géométrie spectrales.
  \par\nopagebreak
  Département de mathématiques, faculté des sciences, université Ibn Tofail, 14000 Kénitra.\par\nopagebreak 
  \textit{E-mail address} : \texttt{nizar.elidrissi@uit.ac.ma}

  \medskip

  \textbf{Hicham Zoubeir.} 
  \par\nopagebreak Laboratoire : Equations aux dérivées partielles, Algèbre et Géométrie spectrales.
  \par\nopagebreak
   Département de mathématiques, faculté des sciences, université Ibn Tofail, 14000 Kénitra.
\par\nopagebreak 
  \textit{E-mail address} : \texttt{hzoubeir2014@gmail.com}
}}

\theoremstyle{plain}
\newtheorem{theorem}{Theorem}[section]
\newtheorem{proposition}{Proposition}[section]

\theoremstyle{definition}

\theoremstyle{remark}
\newtheorem{remark}{Remark}[section]

\makeatletter
\def\keywords{\xdef\@thefnmark{}\@footnotetext}
\makeatother

\begin{document}
\newpage
\maketitle


\begin{abstract}
This note examines the implications of randomly selecting vectors from an infinite-dimensional Hilbert space on linear independence, assuming that for all $k$, the first $k$ vectors follow an absolutely continuous law with respect to a probability measure. It demonstrates that no constraints on the random dimension of their span are necessary, provided that all finite dimensional vector subspaces are considered negligible with respect to the Hilbert space probability measure.
\end{abstract}




\def\thefootnote{*}\footnotetext{These authors contributed equally to this work.}


\tableofcontents

\section{Introduction}

In \cite{ChristensenHasannasab}, Christensen and Hasannasab proved the following result.
\begin{proposition}
\label{PropositionChristensenHasannasab}
Let $E$ be a vector space of (countably or uncountably) infinite dimension and $(e_n)_{n \in \mathbb{N}}$ be a sequence of vectors with values in $E$. Then 
\[ \left( \exists T \in L(\Span\{e_n\}_{n \in \mathbb{N}}, E) : \forall n \in \mathbb{N} : T(e_n)=e_{n+1} \right) \text{ and } \dim\Span\{e_n\}_{n \in \mathbb{N}} = +\infty ] \Rightarrow \]
\[ \quad (e_n)_{n \in \mathbb{N}} \text{ is free.} \] 
\end{proposition}

In \cite{ElIdrissiKabbaj}, the first author of the present article and S. Kabbaj explored a little bit this proposition. They proved for instance that the function appearing in the lower indices of the vectors (namely the successor function $n \in \mathbb{N} \mapsto n+1 \in \mathbb{N})$ can't be replaced by an essentially different function $\varphi : \mathbb{N} \to \mathbb{N}$. They also pushed it in two new directions: they proved that a similar model-theoretic result holds if one accepts to use a certain definition of independence expressible in any algebraic structure, and that the operatorial condition consisting of only one shifting operator can be replaced by a set of operators acting together.

In this short note, we consider a yet different perspective, which is to assume that the vectors are randomly chosen from an infinite dimensional Hilbert space, with no prescription on their mutual law except that for all $k$, the first $k$ vectors are absolutely continuous with respect to the product probability measure. Surprisingly, no condition on the random dimension of their span is needed. However, this comes at the price of assuming that all finite dimensional vector subspaces are negligible with respect to the Hilbert space probability measure. 

\begin{theorem}
\label{MainTheorem}
Let $(\Omega, \mathfrak{F}, P)$ be a probability space and $H$ be a real Hilbert space of infinite dimension. \\
Let $(e_n)_{n \in \mathbb{N}}$ be a sequence of random variables defined on $(\Omega, \mathfrak{F}, P)$ and with values in $H$. \\
Let $Q$ be a complete $\sigma$-finite probability measure on $H$ such that any finite dimensional vector subspace $A$ of $H$ satisfies $Q(A)=0$. \\
Suppose that for all $k \in \mathbb{N}$,  $\left( (e_0,\cdots,e_k)_*(P)  \right)$ is absolutely continuous with respect to $Q^{\otimes (k+1)}$, so that the Radon-Nykodym derivative $\frac{d \left( (e_0,\cdots,e_k)_*(P)  \right)}{dQ^{\otimes (k+1)}}$ exists, which we denote by $p^{(k)} : H^{k+1} \to[0,\infty)$. \\
Then $(e_n)_{n \in \mathbb{N}}$ is free almost surely.
\end{theorem}

For potential probability measures satisfying the above condition, see \cite{VakhaniaTarieladzeChobanyan} and \cite{Kukush}.

\section{Proof of theorem }

\begin{proof}
For all $k \in \mathbb{N}$, let $f^k : (v_0,\cdots,v_k) \in H^{k+1} \mapsto \det \Gram (v_0,\cdots,v_k) \in \mathbb{R}$. \\
By a well-known formula for Gramian determinants exposed in the appendix \ref{FormulaGramianDeterminants}, we have
\[ f^k(v_0,\cdots,v_k) = f^{k-1}(v_0,\cdots,v_{k-1}) h^k(v_0,\cdots,v_{k-1},v_k)^2, \]
where 
\[ h^k(v_0,\cdots,v_k) := d(v_k, \Span(v_0,\cdots,v_{k-1})), \]
and for any $a \in H$ and $V \subseteq H$, $d(a,V)$ denotes the distance of $a$ to $V$. \\
Let's show by induction on $k$ that
\begin{equation}
\label{MeasureReciprocalImagedetGramIsZero}
\forall k \in \mathbb{N} : Q^{\otimes k+1}((f^k)^{-1}(0)) = 0.
\end{equation}
\begin{itemize}
\item If $k=0$, then 
\[ \forall v \in H : f^0(v_0) = \lVert v_0 \rVert^2. \]
Therefore $(f^0)^{-1}(0) = \{0\}$ is $Q$-negligible due to the hypothesis on $Q$, as $\{0\}$ is a finite dimensional vector subspace of $H$.
\item Suppose the property is true at the index $k-1$, $k \geq 1$. \\
Let $v_k \in H$, and suppose that 
\[ \forall (v_0,\cdots,v_{k-1}) \in H^k : f^k(v_0,\cdots,v_{k-1},v_k) = 0. \]
Then 
\[ \forall (v_0,\cdots,v_{k-1}) \in H^k : f^{k-1}(v_0,\cdots,v_{k-1}) h^k(v_0,\cdots,v_k)^2 = 0, \]
and so
\[ \forall (v_0,\cdots,v_{k-1}) \in H^k : \left( h^k(v_0,\cdots,v_{k-1},v_k) = 0 \right) \text{ or } \left( f^{k-1}(v_0,\cdots,v_{k-1}) = 0 \right) \]
Hence 
\[ (f^k)^{-1}(0) \subseteq \left[ (f^{k-1})^{-1}(0) \times H \right] \bigcup \{(v_0,\cdots,v_{k-1},v_k) \in H^{k+1} : h(v_0,\cdots,v_{k-1},v_k) = 0 \}. \]
Therefore, we have
\begin{align*}
0 \leq Q^{\otimes k+1}((f^k)^{-1}(0)) &= \int_{H^{k+1}} \mathbbm{1}_{\{f^k(v_0,\cdots,v_{k-1},v_k)=0\}}  dQ^{\otimes k+1}(v_0,\cdots,v_k) \\
&= \int_{H^{k+1} \setminus \{[(f^{k-1})^{-1}(0) \times H ]\}} \mathbbm{1}_{\{h^k(v_0,\cdots,v_{k-1},v_k)=0\}}  dQ^{\otimes k+1}(v_0,\cdots,v_k) \\
&\quad \textcolor{blue}{\text{since } Q^{\otimes k+1}([(f^{k-1})^{-1}(0) \times H ])=0} \\
&= \int_{H^{k+1}}  \mathbbm{1}_{\{h^k(v_0,\cdots,v_{k-1},v_k)=0\}}  dQ^{\otimes k+1}(v_0,\cdots,v_k) \\
&= \int_{H^k} \left( \int_{H}  \mathbbm{1}_{\{h^k(v_0,\cdots,v_{k-1},v_k)=0\}} dQ(v_k)  \right)  dQ^{\otimes k}(v_0,\cdots,v_{k-1}) \\
&\quad \text{\textcolor{blue}{by Fubini's theorem}} \\
&\leq  \int_{H^k} Q(\Span(v_0,\cdots,v_{k-1}))  dQ^{\otimes k}(v_0,\cdots,v_{k-1}) \\
&= 0
\end{align*}
since $\Span(v_0,\cdots,v_{k-1})$ is a finite dimensional vector subspace of $H$.
\end{itemize}
Hence, equation \ref{MeasureReciprocalImagedetGramIsZero} is proved. \\
Next, let 
\[ \forall k \in \mathbb{N} : E^{(k)} := (e_0, \cdots, e_k), \]
and
\[ \forall k \in \mathbb{N} : G_k := \Gram(e_0,\cdots,e_k) = E^{(k)*}E^{(k)}. \]
By definition of $p^{(k)}$, this function satisfies the property:
\begin{equation}
\label{FormulaExpectedValueg(k)E(k)}
\operatorname E\left[g^{(k)}(E^{(k)})\right] = \int_{H^{k+1}} g^{(k)} (V^{(k)}) p^{(k)}(V^{(k)}) dQ^{\otimes k+1}(V^{(k)})
\end{equation} 
for all bounded Borel measurable $g^{(k)} : H^{k+1} \to \mathbb R$. \\
Notice that we have, for all $\epsilon, t > 0$,
\begin{align*}
P(\det G_k > \epsilon) &= 1 - P(\det G_k \leq \epsilon) \\
&= 1 - P(-\det G_k \geq -\epsilon) \\
&= 1 - P(e^{-t\det G_k} \geq e^{-t\epsilon}) \\
&\geq 1 - \mathbb{E}(e^{-t \det G_k})e^{t\epsilon}
\end{align*}
Passing to the limit $\epsilon \to 0^+$, we obtain
\begin{align*}
P(\det G_k > 0) &\geq 1 - \mathbb{E}(e^{-t \det G_k}) \\
&= 1 - \mathbb{E}(e^{-t \det E^{(k)*}E^{(k)}}) \\
&= 1 - \int_{H^{k+1}} e^{-t \det(V^{(k)*}V^{(k)})} p^{(k)}(V^{(k)}) dQ^{\otimes k+1}(V^{(k)}) \\
&= 1 - \int_{\left[(f^k)^{-1}(0)\right]^c} e^{-t \det(V^{(k)*}V^{(k)})} p^{(k)}(V^{(k)}) dQ^{\otimes k+1}(V^{(k)})
\end{align*}
for any $t > 0$. 
Passing to the limit $t \to +\infty$ and using the dominated convergence theorem, we have
\[ P(\det G_k > 0) \geq 1, \]
hence the equality 
\[ P(\det G_k > 0) = 1. \]
Therefore 
\[ P((e_n)_{n \in \mathbb{N}} \text{ is free}) = P(\det G_k > 0, \text{ for all k } \in \mathbb{N}) = \inf_{k \in \mathbb{N}} P(\det G_k > 0)  = 1. \]
\end{proof}

\begin{remark}
For all $k \in \mathbb{N}$, $\left( (e_0,\cdots,e_k)_*(P)  \right)$ is absolutely continuous with respect to $Q^{\otimes (k+1)}$ in case $(e_n)_{n \in \mathbb{N}}$ is a sequence of independent and absolutely continuous (with respect to $Q$) random variables. Indeed, for any $k \in \mathbb{N}$ and for any measurable subset $V := \bigcup_{j \in \mathbb{N}} V_j := \bigcup_{j \in \mathbb{N}} W^j_0 \times \cdots \times W^j_k \subseteq H^{k+1}$, we have 
\begin{align*}
Q^{\otimes k+1}(V) = 0 &\Rightarrow \forall j \in \mathbb{N} : \exists i \in [\![0,k]\!] : Q(W^j_i) = 0 \\
&\Rightarrow \forall j \in \mathbb{N} : \exists i \in  [\![0,k]\!] : (e_i)_*(P)(W^j_i) = 0 \\
&\Rightarrow  (e_0,\cdots,e_k)_*(P)(V) \leq \sum_{j \in \mathbb{N}} \Pi_{i=0}^k (e_i)_*(P)(W^j_i) = 0,
\end{align*}
which shows that $(e_0,\cdots,e_k)_*(P)$ is absolutely continuous with respect to $Q^{\otimes k+1}$. More precisely, we have $\frac{d \left( (e_0,\cdots,e_u)_*(P)  \right)}{dQ^{\otimes (u+1)}}(v_0,\cdots,v_u) = \Pi_{i=0}^u \frac{d \left((e_i)_*(P)  \right)}{dQ}(v_i)$ in this case.
\end{remark}

\appendix

\section{Ratio of two successive Gramian determinants}
\label{FormulaGramianDeterminants}
Since we couldn't find a good bibliographic reference for this formula, we include its proof in this appendix. \\
Consider $v_0, \cdots, v_k \in H$, and $w_0, \cdots, w_k \in H$ the set produced from the first by iterating the Gram-Schmidt process. \\
For each $m \in [\![0,k]\!]$, the span of the vectors $(w_0,w_1, \cdots,w_m)$ equals the span of the vectors $(v_0,v_1, \cdots, v_m)$.
Let
\[ v'_k := v_k - \sum_{i=0}^{k-1} \langle v_k , w_i \rangle w_i. \]
Since $v'_k$ is orthogonal to all the vectors $w_0, \cdots, w_{k-1}$, it is orthogonal to all the vectors $v_0, \cdots, v_{k-1}$ as well. \\ 
Consider the Gramian determinant
\[ f^k(v_0,\cdots,v_{k-1},v'_k). \]
Using the definition, we see that it is equal to
\[ f^k(v_0,\cdots,v_{k-1},v_k) \]
because $v'_k \in v_k + \Span(v_0,\cdots,v_{k-1})$. \\
But it is also equal to
\[ f^k(v_0,\cdots,v_{k-1}) \lVert v'_k \rVert^2 \]
by a straightforward computation of the determinant of a block-diagonal matrix. \\
Or $\lVert v'_k \rVert$ is nothing else but the explicit formula for $d(v_k, \Span(v_0,\cdots,v_{k-1}))$. \\
Hence, we have the formula:
\[ f^k(v_0,\cdots,v_k) = f^{k-1}(v_0,\cdots,v_{k-1}) h^k(v_0,\cdots,v_{k-1},v_k)^2, \]
where 
\[ h^k(v_0,\cdots,v_k) := d(v_k, \Span(v_0,\cdots,v_{k-1})), \]
and for any $a \in H$ and $V \subseteq H$, $d(a,V)$ denotes the distance of $a$ to $V$.

\section{Finite dimensional affine subspaces are negligible in $H^k$ if they are in $H$}
The result in this appendix is not needed in the proof of the theorem, but since it is a good complement to the subject of our paper and might be of independent interest, we include it here. \\
We claim that under the condition that $Q(A) = 0$ for every finite dimensional affine subspace of $H$, we have $Q^{\otimes k}(A) = 0$ for every strict affine subspace of $H^k$ and for all $k \in \mathbb{N}^*$. \\
Let's show it by induction on $k \in \mathbb{N}^*$:
\begin{itemize}
\item If $k=1$, then $Q(A)=0$ by the hypothesis on $Q$.
\item Suppose the property is true at the index $k$, $k \geq 1$. \\
Let $A$ be a strict affine subspace of $H^{k+1}$. \\
For all $v_k \in H$, let 
\[ C^A(v_k) := \{(v_0,\cdots,v_{k-1}) \in H^k : (v_0,\cdots,v_{k-1},v_k) \in A \}. \]
We have by Fubini's theorem :
\begin{align*}
Q^{\otimes {k+1}}(A) &= \int_{H^{k+1}} \mathbbm{1}_{(v_0,\cdots,v_{k-1},v_k) \in A}  dQ^{\otimes {k+1}}(v_0,\cdots,v_k) \\
&= \int_H \left( \int_{H^k}  \mathbbm{1}_{(v_0,\cdots,v_{k-1}) \in C^A(v_k)}  dQ^{\otimes k}(v_0,\cdots,v_{k-1}) \right) dQ(v_k) \\
&\quad \text{\textcolor{blue}{by Fubini's theorem}} \\
&= \int_{H} Q^{\otimes k}(C^A(v_k)) dQ(v_k) 
\end{align*}
Hence, we see that if the sets $C^A(v_k)$ are proven to be either empty or finite dimensional affine subspaces of $H^k$, we will be done by the induction hypothesis, since we will have $Q^{\otimes k}(C^A(v_k)) = 0$ in the integral. \\
This is the case indeed, since if $C^A(v_k)$ is not empty, then it is the finite dimensional affine subspace passing through $pt \in C^A(v_k)$ and with associated vector subspace $C^A(v_k) - pt$. \\
$C^A(v_k) - pt$ is indeed a finite dimensional vector subspace since
\begin{align*}
x \in C^A(v_k) - pt &\Leftrightarrow pt + x \in C^A(v_k) \\
&\Leftrightarrow (pt + x, v_k) \in A \\
&\Leftrightarrow (pt + x, v_k) \in (pt, v_k) + \vv{A} \\
&\Leftrightarrow (x,0) \in \vv{A}
\end{align*}
where $\vv{A}$ denotes the finite dimensional vector subspace associated to $A$. \\
The last passage explains why we are working in this appendix with affine subspaces (as opposed to vector subspaces), since $C^A(v_k)$ is not in general a vector subspace even if $A$ is. This happens when $A$ doesn't contain $(0,\cdots,0,v_k)$. 
\end{itemize}

\section{Strict affine subspaces are negligible in $H^k$ if they are in $H$}

Again, the result in this appendix is not needed in the proof of the theorem, but since it is a good complement to the subject of our paper and might be of independent interest, we include it here. \\
We claim that under the condition that $Q(A) = 0$ for every strict affine subspace of $H$, we have $Q^{\otimes k}(A) = 0$ for every strict affine subspace of $H^k$ and for all $k \in \mathbb{N}^*$. \\
Let's show it by induction on $k \in \mathbb{N}^*$:
\begin{itemize}
\item If $k=1$, then $Q(A)=0$ by the hypothesis on $Q$.
\item Suppose the property is true at the index $k$, $k \geq 1$. \\
Let $A$ be a strict affine subspace of $H^{k+1}$. \\
For all $v_k \in H$, let 
\[ C^A(v_k) := \{(v_0,\cdots,v_{k-1}) \in H^k : (v_0,\cdots,v_{k-1},v_k) \in A \}. \]
We have by Fubini's theorem :
\begin{align*}
Q^{\otimes {k+1}}(A) &= \int_{H^{k+1}} \mathbbm{1}_{(v_0,\cdots,v_{k-1},v_k) \in A}  dQ^{\otimes {k+1}}(v_0,\cdots,v_k) \\
&= \int_H \left( \int_{H^k}  \mathbbm{1}_{(v_0,\cdots,v_{k-1}) \in C^A(v_k)}  dQ^{\otimes k}(v_0,\cdots,v_{k-1}) \right) dQ(v_k) \\
&\quad \text{\textcolor{blue}{by Fubini's theorem}} \\
&= \int_{H} Q^{\otimes k}(C^A(v_k)) dQ(v_k) 
\end{align*}
Hence, we see that if the sets $C^A(v_k)$ are proven to be either empty or strict affine subspaces of $H^k$ for almost all $v_k$, we will be done by the induction hypothesis, since we will have $Q^{\otimes k}(C^A(v_k)) = 0$ in the integral. \\
First, if $C^A(v_k)$ is not empty, then it is the affine subspace passing through $pt \in C^A(v_k)$ and with associated vector subspace $C^A(v_k) - pt$. \\
$C^A(v_k) - pt$ is indeed a vector subspace since
\begin{align*}
x \in C^A(v_k) - pt &\Leftrightarrow pt + x \in C^A(v_k) \\
&\Leftrightarrow (pt + x, v_k) \in A \\
&\Leftrightarrow (pt + x, v_k) \in (pt, v_k) + \vv{A} \\
&\Leftrightarrow (x,0) \in \vv{A}
\end{align*}
where $\vv{A}$ denotes the vector subspace associated to $A$. \\
The last passage explains why we are working in this appendix with affine subspaces (as opposed to vector subspaces), since $C^A(v_k)$ is not in general a vector subspace even if $A$ is. This happens when $A$ doesn't contain $(0,\cdots,0,v_k)$. \\  
Now, suppose to the contrary that there exists a measurable subset $\alpha \subseteq H$ such that $Q(\alpha)>0$ and $\forall v_k \in \alpha : C^A(v_k) = H^k$. So, $\forall v_k \in \alpha : H^k \times \{v_k\} \subseteq A$. Therefore
\[ H^{k+1} = H^k \times H = H^k \times \langle \{v_k\}_{v_k \in \alpha} \rangle_{\text{aff}} \subseteq A, \]
where $\langle X \rangle_{\text{aff}}$ denotes the smallest affine subspace containing $X$, for any subset $X \subseteq H$, which is a contradiction. \\
Indeed, if $\langle \{v_k\}_{v_k \in \alpha} \rangle_{\text{aff}}$ were a strict affine subspace of $H$, we would have 
\[ 0 \leq Q(\alpha) \leq  Q(\langle \{v_k\}_{v_k \in \alpha} \rangle_{\text{aff}}) = 0, \] 
a contradiction. So, we have $\langle \{v_k\}_{v_k \in \alpha} \rangle_{\text{aff}} = H$ indeed. \\
The last inclusion can be proven as follows. Let $x \in H^k$, $(v_k^j)_{j=1}^l$ in $\alpha$, and $(\alpha_j)_{j=1}^l$ scalars such that $\sum_{j=1}^l \alpha_j = 1$. Then $(x, \sum_{j=1} \alpha_j v_k^j) = \sum_{j=1}^l \alpha_j (x,v_k^j) \in A$, since $A$ is affine.
\end{itemize}

\nocite{*}
\bibliographystyle{alpha}
\bibliography{references}

\Addresses

\end{document}